\documentclass[11pt, a4paper,pagebackref]{amsart}
\usepackage{amsmath}
\usepackage{amssymb}
\usepackage{amsthm}
\usepackage{amscd}
\usepackage{tikz-cd}
\usepackage{color}
\usepackage{hyperref}

\usepackage[inner=2.3cm,outer=2.3cm, bottom=3.3cm]{geometry}
\usepackage{setspace}

\numberwithin{equation}{section}

\usepackage[colorinlistoftodos]{todonotes}

\usepackage[all]{xy}
\usepackage{enumerate}
\usepackage{hyperref}

\usepackage{color}
\hypersetup{colorlinks=true} 

\newtheorem{theorem}{Theorem}[section]

\newtheorem{proposition}[theorem]{Proposition}
\newtheorem{corollary}[theorem]{Corollary}

\theoremstyle{definition}
\newtheorem{definition}[theorem]{Definition}
\newtheorem{newclaim}[theorem]{}

\newtheorem{remark}[theorem]{Remark}
\newtheorem{remark and definition}[theorem]{Remark and Definition}
\newtheorem{remark and notation}[theorem]{Remark and Notation}

\newtheorem*{fund}{Funding}
\newtheorem*{agra}{Acknowlegment}

\newtheorem{fact}[theorem]{Fact}

\newcommand\Spec{\operatorname{Spec}}

\newcommand\Tor{\operatorname{Tor}}

\newcommand\depth{\operatorname{depth}}

\newcommand\Ker{\operatorname{\Ker}}

\newcommand\id{\operatorname{id}}
\newcommand\Supp{\operatorname{Supp}}

\newcommand{\ldt}{\otimes^\mathbf{L}}

\author[Martins]{Paulo Martins}

\address{Universidade de S{\~a}o Paulo -
ICMC, Caixa Postal 668, 13560-970, S{\~a}o Carlos-SP, Brazil}
\email{paulomartinsmtm@gmail.com}

\keywords{Chouinard's formula, Bass' formula, $C$-quasi-injective dimension, $C$-injective dimension, quasi-injective dimension, semidualizing modules}
\subjclass[2020]{13D05, 13D07.}

\begin{document}

\title{Chouinard's formula for $C$-quasi-injective dimension}

\maketitle

\begin{abstract}

The $C$-quasi-injective dimension is a recently introduced homological invariant that unifies and extends the notions of quasi-injective dimension and of injective dimension with respect to a semidualizing module, previously studied by Gheibi and by Takahashi and White, respectively. In the main results of this paper, we provide extensions of the Bass' formula and a version of the Chouinard's formula for modules of finite $C$-quasi-injective dimension over an arbitratry ring. \end{abstract}
\section{Introduction}
Throughout this paper, all rings are assumed to be commutative and Noetherian. In 1976, Chouinard provided a general formula for the injective dimension of a module, whenever it is finite (see \cite{Chouinard}), without assuming that the base ring is local or that the module is finitely generated. Let $M$ be an $R$-module with finite injective dimension, Chouinard's formula states that 
\begin{align*}
    \operatorname{id}_R M = \sup \lbrace \operatorname{depth} R_{\mathfrak{p}}- \operatorname{width}_{R_{\mathfrak{p}}} M_{\mathfrak{p}} \mid \mathfrak{p} \in \operatorname{Spec}R \rbrace.
\end{align*}
Later, Khatami, Tousi and Yassemi \cite{yassemi} proved a version of Chouinard's formula for Gorenstein-injective dimension. We recall that for an $R$-module $M$ over a local ring $R$ with residue field $k$, $\operatorname{width}_R M$ is equal to $\inf \lbrace i \mid \operatorname{Tor}_i^R(k,M) \neq 0 \rbrace$.

The quasi-injective dimension (qid) is a refinement of the classical notion of  injective dimension of a module, in the sense that there is always an inequality $\operatorname{qid}_R M \leq \operatorname{id}_R M$ and equality holds when $\operatorname{id}_RM$ is finite. It was introduced by Gheibi \cite{Gheibi} that recovered several well-known results about injective and Gorenstein-injective dimensions in the context of quasi-injective dimension. In a recent paper, Tri \cite{Tri} obtained a Chouinard's formula for quasi-injective dimension and then extended the Bass' Formula for quasi-injective dimension previously proved by Gheibi \cite[Theorem 3.2]{Gheibi}. Moreover, Jorge-Pérez, Martins and Mendoza-Rubio proved a version of Ischebeck's formula for quasi-injective dimension (see \cite[Theorem 4.3]{VH-P-VD}) .

Let $C$ be a semidualizing $R$-module. Recently, Dey, Ferraro, and Gheibi \cite{ferraro} introduced the notion of $C$-quasi-injective dimension, unifying and extending both the theory of $C$-injective dimension developed by Takahashi and White \cite{takahashi} and the theory of quasi-injective dimension introduced by Gheibi \cite{Gheibi}. More precisely, for every $R$-module $M$, one has $\operatorname{\textit{C}-qid}_R M \leq \operatorname{\textit{C}-id}_R M$, and, in the special case $C=R$, the notion of $C$-quasi-injective dimension recovers the quasi-injective dimension defined in \cite{Gheibi}.

It was shown in \cite[Theorem 7.4]{ferraro} that the Bass formula for $C$-quasi-injective dimension holds for finitely generated modules of finite $C$-quasi-injective dimension over a local ring $R$ whenever either $R$ is Cohen-Macaulay or $\operatorname{Tor}_{>0}^R(C,M)=0$. Furthermore, a version of Ischebeck's formula was established in \cite[Theorem 7.7]{ferraro} under the assumptions that the pair of finitely generated $R$-modules $(M,N)$ satisfies $\operatorname{Ext}_R^{\gg 0}(M,N)=0$, $M$ belongs to the Auslander class $\mathcal{A}_C(R)$, $\operatorname{Tor}_{>0}^R(C,N)=0$, and $\operatorname{\textit{C}-qid}_R N<\infty$.

Motivated by the articles previously mentioned, the objective of this paper is to provide a Chouinard's formula for $C$-quasi-injective dimension. The main results of this paper extend and unify the Bass formula for $C$-quasi-injective dimension proved in \cite[Theorem 7.4]{ferraro} when $\operatorname{Tor}_{>0}^R(C,M)=0$ and the recently proved Chouinard formula for quasi-injective dimension \cite[Theorem 3.3]{Tri}. Let us now briefly describe the contents of this paper. In Section \ref{preli}, we begin by introducing the necessary notations, definitions, and foundational results that will be utilized throughout the paper. In Section \ref{main}, we first establish a version of Chouinard's formula for $C$-injective dimension (Proposition \ref{lemma1}) and then prove the following main results.
\begin{theorem}[See Theorem \ref{Cchouinard}]\label{teo1.1}
Let $C$ be a semidualizing $R$-module and let $M$ be an $R$-module of finite and positive $C$-quasi-injective dimension with $\operatorname{Tor}_{>0}^R(C,M)=0$. Then
$$\operatorname{\textit{C}-qid}_RM = \sup \lbrace \operatorname{depth} R_{\mathfrak{p}}- \operatorname{width}_{R_{\mathfrak{p}}} M_{\mathfrak{p}} \mid \mathfrak{p} \in \operatorname{Spec}R \rbrace.$$
\end{theorem}
\begin{theorem}[See Theorem \ref{theoremfg}]
    Let $C$ be a semidualizing $R$-module and let $M$ be a finitely generated $R$-module of finite $C$-quasi-injective dimension with $\operatorname{Tor}_{>0}^R(C,M)=0$. Then
$$\operatorname{\textit{C}-qid}_RM = \sup \lbrace \operatorname{depth} R_{\mathfrak{p}} \mid \mathfrak{p} \in \operatorname{Supp}_RM \rbrace.$$
\end{theorem}
As a corollary of Theorem \ref{teo1.1}, we show that the $C$-quasi-injective dimension of a module coincides with the $C$-injective dimension whenever the latter is finite (see Corollary \ref{crl:comparativo}). In the final section, we obtain a criterion for finiteness of $C$-injective dimension that is a dual version of \cite[Theorem 6.11]{ferraro} and recovers \cite[Theorem 4.6]{Gheibi} when $C=R$.
\begin{corollary}[See Corollary \ref{crl}]
Let $C$ be a semidualizing $R$-module. If $M$ is an $R$-module such that
\begin{enumerate}
    \item $\operatorname{\textit{C}-qid}_R M< \infty$,
    \item $M \in \mathcal{A}_C(R)$,
    \item $\operatorname{Ext}_R^{>0}(M,M)=0$,
\end{enumerate}
then $\operatorname{\textit{C}-id}_R M< \infty$.
\end{corollary}
\section{Preliminaries}\label{preli}
In this section, we introduce fundamental definitions and facts that will be considered throughout the paper. All complexes in this paper are indexed homologically. 
\begin{newclaim}    
For a complex $$X_{\bullet}=( \cdots \stackrel{\partial_{i+2}}{\longrightarrow} X_{i+1} \stackrel{\partial_{i+1}}{\longrightarrow  } X_i \stackrel{\partial_{i}}{\longrightarrow
      } X_{i-1}   \longrightarrow  \cdots) $$ of $R$-modules, we set for each integer $i$, $Z_i(X_{\bullet})=\ker \partial_i$ and $B_i(X_{\bullet})= \operatorname{Im} \partial_{i+1}$ and $H_i(X_{\bullet})=Z_i(X_{\bullet})/B_i(X_{\bullet})$.  Moreover, we set: \begin{align*}
&\left\{ 
    \begin{aligned}
        \sup X_{\bullet} &= \sup \{ i \in \mathbb{Z} : X_i \neq 0 \},\\
        \inf X_{\bullet} &= \inf \{ i \in \mathbb{Z} : X_i \neq 0 \},
    \end{aligned}
\right. \quad 
\left\{
    \begin{aligned}
        \operatorname{hsup } X_{\bullet} &= \sup \{ i \in \mathbb{Z} : H_i(X_{\bullet}) \neq 0 \},\\
        \operatorname{hinf } X_{\bullet} &= \inf \{ i \in \mathbb{Z} : H_i(X_{\bullet}) \neq 0 \}.
    \end{aligned}
\right.
\end{align*}
The \textit{length} of $X_{\bullet}$ is defined to be  $\operatorname{length} X_{\bullet}= \sup X_{\bullet} - \inf X_{\bullet}$. We say that $X_{\bullet}$ is \textit{bounded}, if $\operatorname{length} X_{\bullet} < \infty$. We say that $X_{\bullet}$ is \textit{bounded above} if $\sup X_{\bullet} < \infty$.
\end{newclaim}
\begin{newclaim}[\textbf{Small restricted injective dimension and width of modules}] Let $M$ be an $R$-module. The \textit{small restricted injective dimension}  of $M$, denoted by $\operatorname{rid}_R M$, is defined as follows
\begin{align*}
    \operatorname{rid}_R M = \sup \lbrace i \in \mathbb{N}_0 \mid & \operatorname{Ext}_R^i(N,M) \neq 0 \text{ for some finitely generated} \\ & R\text{-module } N \text{ of finite projective dimension} \rbrace.
\end{align*}
Let $\mathfrak{a}$ be an ideal of $R$ generated by $\boldsymbol{a}=a_1,\dots,a_t$. The \textit{$\mathfrak{a}$-width} of $M$ is defined as follows: 
\begin{align*}
\operatorname{width}_R (\mathfrak{a},M) = \inf \lbrace i \in \mathbb{N}_0 \mid H_i(M \otimes_R K(\boldsymbol{a})) \neq 0 \rbrace
\end{align*}
where $K(\boldsymbol{a})$ is the Koszul complex. Moreover, when $R$ is a local ring with maximal ideal $\mathfrak{m}$, we set $\operatorname{width}_R M: = \operatorname{width}_R (\mathfrak{m},M)$.
\end{newclaim}
\begin{fact}\cite[Proposition 4.9]{lars2} Let $\mathfrak{a}$ be an ideal of $R$ and let $M$ be an $R$-module. Then
    \begin{align*}
        \operatorname{width}_R (\mathfrak{a},M) = \inf \lbrace i \in \mathbb{N}_0 \mid \operatorname{Tor}_i^R(R/\mathfrak{a},M) \neq 0 \rbrace.
    \end{align*}
\end{fact}
It is easy to see that if $M$ is a finitely generated module over a local ring, then $\operatorname{width}_R(\mathfrak{a},M)=0$ for any non-zero ideal $\mathfrak{a}$. In particular, $\operatorname{width}_R M=0$. 

The following fact follows directly by \cite[Proposition 5.3(c)]{lars2}.

\begin{fact}\label{factrid}
Let $R$ be a local ring and let $M$ be an $R$-module. If $\operatorname{rid}_R M \leq 0$, then $$\operatorname{depth}R- \operatorname{width}_R M \leq 0.$$
\end{fact}

We refer the reader to \cite{lars2} for details about small restricted injective dimension and width.
\begin{newclaim}[\textbf{Semidualizing modules}] A finitely generated $R$-module $C$ is called a \textit{semidualizing} $R$-module if 
\begin{enumerate}
    \item The natural homothety map $R \to \operatorname{Hom}_R(C,C)$ is an isomorphism.
    \item $\operatorname{Ext}_R^i(C,C)=0$ for all $i>0$.
\end{enumerate}
\end{newclaim}
\begin{newclaim}[\textbf{Auslander and Bass classes}]
Let $C$ be a semidualizing $R$-module.
The \textit{Auslander class} $\mathcal{A}_C (R)$ is the class of $R$-modules $M$ satisfying in the following conditions.
\begin{enumerate}
    \item The natural map $M \rightarrow \operatorname{Hom}_R(C, C \otimes_R M)$ is an isomorphism. 
    \item One has $\operatorname{Tor}_{>0}^R(C,M)=0=\operatorname{Ext}_R^{>0}(C,C \otimes_R M)$.
\end{enumerate}
The \textit{Bass class} $\mathcal{B}_C(R)$ is the class of $R$-modules $M$ satisfying in the following conditions.
\begin{enumerate}
    \item The evaluation map $C \otimes_R \operatorname{Hom}_R(C,M) \rightarrow M$ is an isomorphism.
    \item One has $\operatorname{Ext}_{R}^{>0} (C,M)=0=\operatorname{Tor}_{>0}^R(C,\operatorname{Hom}_R(C,M))$.
\end{enumerate}
\end{newclaim}
\begin{newclaim}[\textbf{$C$-injective dimension}] Let $C$ be a semidualizing $R$-module and let $I$ be an injective $R$-module. The module $\operatorname{Hom}_R (C,I)$ is called a \textit{$C$-injective} $R$-module. For an $R$-module $M$, a \textit{$C$-injective resolution} of $M$ is an exact complex
\begin{align*}
   0 \rightarrow M \rightarrow \operatorname{Hom}_R(C,I_0) \rightarrow \operatorname{Hom}_R(C,I_1) \rightarrow \cdots, 
\end{align*}
where the $I_i$'s are injective $R$-modules. We say $\operatorname{\textit{C}-id}_R M<\infty$ if $M$ admits a bounded $C$-injective resolution. Moreover, we say that $\operatorname{\textit{C}-id}_R M =n$ if the smallest $C$-injective resolution of $M$ has length $n$. 
\end{newclaim}

\begin{theorem}\cite[Corollary 2.9(b)]{takahashi} \label{corollary2.9}
Let $C$ be a semidualizing $R$-module. The class $\mathcal{A}_C(R)$ contains  every $R$-module of finite $C$-injective dimension
\end{theorem}

\begin{theorem}\cite[Theorem 2.11(b)]{takahashi}\label{theorem2.11}
Let $C$ be a semidualizing $R$-module and let $M$ be an $R$-module. Then $\operatorname{\textit{C}-id}_R M= \id_R(C \otimes_RM)$.
\end{theorem}

The author refer to the classical references \cite{takahashi,SemidualizingModules} for more about semidualizing modules, Auslander and Bass classes and $C$-injective dimension. 

The definition of $C$-quasi-injective dimension was recently introduced by Dey, Ferraro and Gheibi \cite{ferraro}. This definition extend and unify the theories of $C$-injective and quasi-injective dimensions.

\begin{definition}
    Let $C$ be a semidualizing $R$-module. An $R$-module $M$ is said to have finite $C$-quasi-injective dimension if there exists a bounded complex $I_{\bullet}$ of injective $R$-modules such that $\operatorname{Hom}_R(C,I_{\bullet})$ is not aclyclic and all the homologies are finite direct sum of copies of $M$ (or zero). Such a complex $I_{\bullet}$ is said to be a \textit{$C$-quasi-injective resolution} of $M$. The \textit{$C$-quasi-injective dimension} of $M$ is defined as
     \begin{align*}
            \operatorname{\textit{C}-qid}_R M = \inf \lbrace \operatorname{hinf} (\operatorname{Hom}_R(C,I_{\bullet})) - \operatorname{inf} (\operatorname{Hom}_R(C,I_{\bullet})) :  I_{\bullet} \text{ is a $C$-quasi-injective resolution of } M \},
        \end{align*}
        if  $M\not=0$, and $\operatorname{\textit{C}-qid}_R M=-\infty$ if $M=0$.
\end{definition}
We remark that when $C=R$, the above definition recovers the quasi-injective dimension introduced by Gheibi in \cite{Gheibi}. Also, it is easy to see that every module $M$ of finite $C$-injective dimension has finite $C$-quasi-injective and that $\operatorname{\textit{C}-qid}_RM \leq \operatorname{\textit{C}-id}_R M$. We show in Corollary \ref{crl:comparativo} that equality holds provided that $\operatorname{\textit{C}-id}_R M$ is finite.
\section{Main results}\label{main}
In this section, we prove the main result of the paper by establishing a Chouinard's formula for the $C$-quasi-injective dimension. In our proofs, all complexes are indexed homologically.

We begin with the following proposition, where we establish a Chouinard's formula for modules of finite injective dimension relative to a semidualizing module. This result will play a key role in the proof of the main theorem. Recall that if $M$ is an $R$-module over a local ring $R$ with residue field $k$, then
\(
\operatorname{width}_R M = \inf \left\{ i \in \mathbb{N}_0 \mid \operatorname{Tor}_i^R(k,M) \neq 0 \right\}.
\)
\begin{proposition}\label{lemma1}
Let $C$ be a semidualizing $R$-module and let $M$ be an $R$-module of finite $C$-injective dimension. Then 
$$\operatorname{\textit{C}-id}_RM = \sup \lbrace \operatorname{depth} R_{\mathfrak{p}}- \operatorname{width}_{R_{\mathfrak{p}}} M_{\mathfrak{p}} \mid \mathfrak{p} \in \operatorname{Spec}R \rbrace.$$
\begin{proof}
Since $\operatorname{\textit{C}-id}_RM< \infty$, $M \in \mathcal{A}_C(R)$ by Theorem \ref{corollary2.9}. Then, we have  $\Tor_{>0}^R(C,M)=0$. By Theorem \ref{theorem2.11}, we have $\operatorname{\textit{C}-id}_RM= \operatorname{id}_R (C \otimes_R M)$. So, by the classical Chouinard's formula \cite[Corollary 3.1]{Chouinard}, we have
\begin{align*}
\operatorname{\textit{C}-id}_RM & = \operatorname{id}_R (C \otimes_R M) \\
& = \sup \lbrace \operatorname{depth} R_{\mathfrak{p}}- \operatorname{width}_{R_{\mathfrak{p}}} ( C_{\mathfrak{p}} \otimes_{R_{\mathfrak{p}}} M_{\mathfrak{p}}) \mid \mathfrak{p} \in \operatorname{Spec}R \rbrace.\end{align*}
To finish this proof, we prove that $\operatorname{width}_{R_{\mathfrak{p}}} ( C_{\mathfrak{p}} \otimes_{R_{\mathfrak{p}}} M_{\mathfrak{p}}) = \operatorname{width} M_{\mathfrak{p}}$ for all $\mathfrak{p} \in \Spec R$. For each $\mathfrak{p} \in \Spec R$, $\operatorname{Tor}_{>0}^{R_{\mathfrak{p}}}(C_{\mathfrak{p}},M_{\mathfrak{p}})=0$ and therefore $C_{\mathfrak{p}}\otimes_{R_{\mathfrak{p}}} M_{\mathfrak{p}} $ and $ C_{\mathfrak{p}} \ldt_{R_{\mathfrak{p}}} M_{\mathfrak{p}}$ are isomorphic in the derived category of $R_{\mathfrak{p}}$, by \cite[Proposition 7.4.22]{lars}.
It then follows from \cite[Theorem 16.2.9]{lars} that as $C_{\mathfrak{p}}$ is a finitely generated $R_{\mathfrak{p}}$-module (i.e, $\operatorname{width}_{R_{\mathfrak{p}} }C_{\mathfrak{p}}=0$), $\operatorname{width}_{R_{\mathfrak{p}}}  M_{\mathfrak{p}} $ and $\operatorname{width}_{R_{\mathfrak{p}}} ( C_{\mathfrak{p}} \otimes_{R_{\mathfrak{p}}} M_{\mathfrak{p}}) $ are simultaneously finite and $\operatorname{width}_{R_{\mathfrak{p}}} ( C_{\mathfrak{p}} \otimes_{R_{\mathfrak{p}}} M_{\mathfrak{p}}) =\operatorname{width}_{R_{\mathfrak{p}}}  M_{\mathfrak{p}} $ for all $\mathfrak{p} \in \Spec R$. Therefore, the desired equality follows.
\end{proof}
\end{proposition}
Now, we can prove the main theorem of this paper, which extends Chouinard’s Formula to the setting of $C$-quasi-injective dimension. We remark that the assumption $\Tor_{>0}^R(C,M)=0$ seems to be natural when working with $C$-quasi-injective dimension, as can be seen, for instance, in \cite[Proposition 7.3, Theorem 7.4 and Theorem 7.7]{ferraro}.
\begin{theorem}\label{Cchouinard}
Let $C$ be a semidualizing $R$-module and let $M$ be an $R$-module of finite and positive $C$-quasi-injective dimension  with $\operatorname{Tor}_{>0}^R(C,M)=0$. Then
$$\operatorname{\textit{C}-qid}_RM = \sup \lbrace \operatorname{depth} R_{\mathfrak{p}}- \operatorname{width}_{R_{\mathfrak{p}}} M_{\mathfrak{p}} \mid \mathfrak{p} \in \operatorname{Spec}R \rbrace.$$
\begin{proof}
Let $I_{\bullet}$ be a bounded $C$-quasi-injective resolution of $M$ with $\operatorname{\textit{C}-qid}_RM = \operatorname{hinf}(\operatorname{Hom}_R(C,I_{\bullet})) - \operatorname{inf}(\operatorname{Hom}_R(C,I_{\bullet}))$. Set $s= \operatorname{hinf}(\operatorname{Hom}_R(C,I_{\bullet}))$. There are exact sequences
\begin{align}\label{nonlocal}
\begin{cases}
0 \rightarrow Z_i \rightarrow \operatorname{Hom}_R(C,I_i) \rightarrow B_{i-1} \rightarrow 0 \\  0 \rightarrow B_i \rightarrow Z_i \rightarrow H_i(\operatorname{Hom}_R(C,I_{\bullet})) \rightarrow 0
\end{cases}
\quad (i \in \mathbb{Z})
\end{align}
where $Z_i=\operatorname{ker}(\partial_i^{\operatorname{Hom}_R(C,I_{\bullet})})$, $B_i=\operatorname{Im}(\partial_{i+1}^{\operatorname{Hom}_R(C,I_{\bullet})})$ and $H_i(\operatorname{Hom}_R(C,I_{\bullet})) \cong M^{\oplus b_i}$ for some $b_i \geq 0$. It is clear that $\operatorname{\textit{C}-qid}_RM=\operatorname{\textit{C}-id}_R Z_s$.  For any $\mathfrak{p} \in \operatorname{Spec} R$, we have the exact sequences of $R_{\mathfrak{p}}$-modules
\begin{align}\label{localized}
\begin{cases}
0 \rightarrow (Z_i)_{\mathfrak{p}} \rightarrow \operatorname{Hom}_R(C,I_i)_{\mathfrak{p}} \rightarrow (B_{i-1})_{\mathfrak{p}} \rightarrow 0 \\  0 \rightarrow (B_i)_{\mathfrak{p}} \rightarrow (Z_i)_{\mathfrak{p}} \rightarrow M_{\mathfrak{p}}^{\oplus b_i} \rightarrow 0
\end{cases}
\quad (i \in \mathbb{Z}).
\end{align}
The exact sequence $0 \rightarrow (B_s)_{\mathfrak{p}} \rightarrow (Z_s)_{\mathfrak{p}} \rightarrow M_{\mathfrak{p}}^{\oplus b_s} \rightarrow 0 $ ($b_s>0$) induces the long exact sequence
\begin{align*}
   \cdots \rightarrow \operatorname{Tor}_i^{R_{\mathfrak{p}}}(k(\mathfrak{p}),(B_s)_{\mathfrak{p}}) \rightarrow \operatorname{Tor}_i^{R_{\mathfrak{p}}}(k(\mathfrak{p}),(Z_s)_{\mathfrak{p}}) \\
   \rightarrow \operatorname{Tor}_i^{R_{\mathfrak{p}}} (k(\mathfrak{p}),M_{\mathfrak{p}})^{\oplus b_s} & \rightarrow \operatorname{Tor}_{i-1}^{R_{\mathfrak{p}}} (k(\mathfrak{p}),(B_s)_{\mathfrak{p}}) \rightarrow \cdots
\end{align*}
where $k(\mathfrak{p})$ denotes the residue field of $R_{\mathfrak{p}}$. Therefore, we have the following inequalities:
\begin{align}
\operatorname{width}_{R_{\mathfrak{p}}} ((Z_s)_{\mathfrak{p}}) & \geq \min \lbrace \operatorname{width}_{R_{\mathfrak{p}}} M_{\mathfrak{p}}, \operatorname{width}_{R_{\mathfrak{p}}} ((B_s)_{\mathfrak{p}}) \rbrace \label{ineq3} \\
\operatorname{width}_{R_{\mathfrak{p}}} M_{\mathfrak{p}} & \geq \min \lbrace \operatorname{width}_{R_{\mathfrak{p}}} ((Z_s)_{\mathfrak{p}}), \operatorname{width}_{R_{\mathfrak{p}}} ((B_s)_{\mathfrak{p}})+1 \rbrace \label{ineq4}
\end{align}

Suppose that $\mathfrak{q} \in \operatorname{Spec} R$ is such that $\operatorname{width}_{R_{\mathfrak{q}}} ((B_s)_{\mathfrak{q}}) \leq \operatorname{width}_{R_{\mathfrak{q}}} M_{\mathfrak{q}}$. By (\ref{ineq3}), we see that $\operatorname{width}_{R_{\mathfrak{q}}} ((B_s)_{\mathfrak{q}}) \leq \operatorname{width}_{R_{\mathfrak{q}}}((Z_s)_{\mathfrak{q}})$.

\text{We claim that:} If $\mathfrak{q} \in \operatorname{Spec} R$ and $\operatorname{width}_{R_{\mathfrak{q}}} ((B_s)_{\mathfrak{q}}) \leq \operatorname{width}_{R_{\mathfrak{q}}} M_{\mathfrak{q}}$, then $\operatorname{depth} R_{\mathfrak{q}} \leq \operatorname{width}_{R_{\mathfrak{q}}} ((B_s)_{\mathfrak{q}})$. 

\textit{Proof of Claim:} By contradiction, we assume that $$d:=\operatorname{depth} R_{\mathfrak{q}}> \operatorname{width}_{R_{\mathfrak{q}}} ((B_s)_{\mathfrak{q}})=:w.$$ Since $I_{s+1}$ is injective, $\operatorname{Hom}_R(C,I_{s+1})$ is $C$-injective and therefore Proposition \ref{lemma1} gives the inequality $\operatorname{width}_{R_{\mathfrak{q}}} (\operatorname{Hom}_R(C,I_{s+1})_{\mathfrak{q}} ) \geq d$. That is, we have $\operatorname{Tor}_i^{R_{\mathfrak{q}}}(k(\mathfrak{q}),\operatorname{Hom}_R(C,I_{s+1})_{\mathfrak{q}})=0$ for $i=0,1,\dots,d-1$. Assume $w>0$. Therefore the exact sequence 
$$0 \rightarrow (Z_{s+1})_{\mathfrak{q}} \rightarrow \operatorname{Hom}_R(C,I_{s+1})_{\mathfrak{q}} \rightarrow (B_{s})_{\mathfrak{q}} \rightarrow 0$$
induces the following exact sequence:
\begin{align*}
\operatorname{Tor}_w^{R_{\mathfrak{q}}}(k(\mathfrak{q}),\operatorname{Hom}_R(C,I_{s+1})_{\mathfrak{q}})=0 \rightarrow \operatorname{Tor}_w^{R_{\mathfrak{q}}}(k(\mathfrak{q}),(B_s)_{\mathfrak{q}}) \rightarrow \operatorname{Tor}_{w-1}^{R_{\mathfrak{q}}}(k(\mathfrak{q}),(Z_{s+1})_{\mathfrak{q}}).
\end{align*}
Then, $\operatorname{Tor}_{w-1}^{R_{\mathfrak{q}}}(k(\mathfrak{q}),(Z_{s+1})_{\mathfrak{q}}) \neq0$, as $\operatorname{Tor}_w^{R_{\mathfrak{q}}}(k(\mathfrak{q}),(B_s)_{\mathfrak{q}}) \neq 0$. Since we are assuming that $ w=\operatorname{width}_{R_{\mathfrak{q}}} ((B_s)_{\mathfrak{q}}) \leq \operatorname{width}_{R_{\mathfrak{q}}} M_{\mathfrak{q}}$, it follows that $\operatorname{Tor}_i^{R_{\mathfrak{q}}} (k(\mathfrak{q}),M_{\mathfrak{q}})=0$ for $i=0,1,\dots, w-1$, and the exact sequence $0 \rightarrow (B_{s+1})_{\mathfrak{q}} \rightarrow (Z_{s+1})_{\mathfrak{q}} \rightarrow M_{\mathfrak{q}}^{\oplus b_{s+1}} \rightarrow 0$ induces the following exact sequence
$$\operatorname{Tor}_{w-1}^{R_{\mathfrak{q}}}(k(\mathfrak{q}),(B_{s+1})_{\mathfrak{q}}) \rightarrow \operatorname{Tor}_{w-1}^{R_{\mathfrak{q}}}(k(\mathfrak{q}),(Z_{s+1})_{\mathfrak{q}}) \rightarrow0$$
and then $\operatorname{Tor}_{w-1}^{R_{\mathfrak{q}}}(k(\mathfrak{q}),(B_{s+1})_{\mathfrak{q}}) \neq 0$, since $\operatorname{Tor}_{w-1}^{R_{\mathfrak{q}}}(k(\mathfrak{q}),(Z_{s+1})_{\mathfrak{q}}) \neq 0$. Repeating this argument a finite number of times using the localized short exact sequences (\ref{localized}), we obtain that $k(\mathfrak{q}) \otimes_{R_{\mathfrak{q}}} (B_{s+w})_{\mathfrak{q}} \neq 0 .$ Moreover, by a previous argument we see that $\operatorname{width}_{R_{\mathfrak{q}}} (\operatorname{Hom}_{R}(C,I_{s+w+1})_{\mathfrak{q}})\geq d$. Then $k(\mathfrak{q}) \otimes_{R_{\mathfrak{q}}} \operatorname{Hom}_R(C,I_{s+w+1})_{\mathfrak{q}} =0$. Therefore, applying $k(\mathfrak{q}) \otimes_{R_{\mathfrak{q}}}-$ to the exact sequence $$0 \rightarrow (Z_{s+w+1})_{\mathfrak{q}} \rightarrow \operatorname{Hom}_{R}(C,I_{s+w+1})_{\mathfrak{q}} \rightarrow (B_{s+w})_{\mathfrak{q}} \rightarrow 0$$ we have a contradiction, since $k(\mathfrak{q}) \otimes_{R_{\mathfrak{q}}} (B_{s+w})_{\mathfrak{q}} \neq 0$. Note that, if $w=0$, then we have the same contradiction by simply applying $k(\mathfrak{q}) \otimes_{R_{\mathfrak{q}}}-$ to the exact sequence $$0 \rightarrow (Z_{s+1})_{\mathfrak{q}} \rightarrow \operatorname{Hom}_R(C,I_{s+1})_{\mathfrak{q}} \rightarrow (B_{s})_{\mathfrak{q}} \rightarrow 0.$$ Thus, $\operatorname{depth} R_{\mathfrak{q}} \leq \operatorname{width}_{R_{\mathfrak{q}}} ((B_s)_{\mathfrak{q}})$ for each $\mathfrak{q} \in \operatorname{Spec} R$ such that $\operatorname{width}_{R_{\mathfrak{q}}} ((B_s)_{\mathfrak{q}}) \leq \operatorname{width}_{R_{\mathfrak{q}}} M_{\mathfrak{q}}$ and the Claim is proved. 

\noindent Therefore, each $\mathfrak{q} \in \Spec R$ such that $\operatorname{width}_{R_{\mathfrak{q}}} ((B_s)_{\mathfrak{q}}) \leq \operatorname{width}_{R_{\mathfrak{q}}} M_{\mathfrak{q}}$ satisfies the following inequalities:
\begin{align}\label{ineq2}
   \operatorname{depth} R_{\mathfrak{q}}-  \operatorname{width}_{R_{\mathfrak{q}}} M_{\mathfrak{q}} & \leq 0 \\
\label{ineq1}
\operatorname{depth} R_{\mathfrak{q}} - \operatorname{width}_{R_{\mathfrak{q}}} ((Z_s)_{\mathfrak{q}}) & \leq 0 .
\end{align}
Using the inequality (\ref{ineq1}) and Proposition \ref{lemma1}, since $\operatorname{\textit{C}-id}_{R}Z_s = \operatorname{\textit{C}-qid}_R M>0$, we have 
\begin{align*}
    \operatorname{\textit{C}-qid}_R M &  = \operatorname{\textit{C}-id}_R Z_s  \\
    & = \sup \lbrace \operatorname{depth} R_{\mathfrak{p}} - \operatorname{width}_{R_{\mathfrak{q}}} ((Z_s)_{\mathfrak{p}}) \mid \mathfrak{p} \in \operatorname{Spec} R \rbrace \\
    & = \sup \lbrace \operatorname{depth} R_{\mathfrak{p}} - \operatorname{width}_{R_{\mathfrak{q}}} ((Z_s)_{\mathfrak{p}}) \mid \mathfrak{p} \text{ with } \operatorname{width}_{R_{\mathfrak{p}}} ((B_s)_{\mathfrak{p}}) > \operatorname{width}_{R_{\mathfrak{p}}} M_{\mathfrak{p}} \rbrace.
\end{align*}
Using (\ref{ineq3}) and (\ref{ineq4}), we see that for any $\mathfrak{p} \in \operatorname{Spec} R$ satisfying $\operatorname{width}_{R_{\mathfrak{p}}} ((B_s)_{\mathfrak{p}}) > \operatorname{width}_{R_{\mathfrak{p}}} M_{\mathfrak{p}}$, it follows that $\operatorname{width}_{R_{\mathfrak{p}}} M_{\mathfrak{p}}=\operatorname{width}_{R_{\mathfrak{p}}} ((Z_s)_{\mathfrak{p}})$. Finally, since $\operatorname{\textit{C}-qid}_RM>0$, using  (\ref{ineq2}) we have
\begin{align*}
\operatorname{\textit{C}-qid}_R M & = \sup \lbrace \operatorname{depth} R_{\mathfrak{p}} - \operatorname{width}_{R_{\mathfrak{q}}} ((Z_s)_{\mathfrak{p}}) \mid \mathfrak{p} \text{ with } \operatorname{width}_{R_{\mathfrak{p}}} ((B_s)_{\mathfrak{p}}) > \operatorname{width}_{R_{\mathfrak{p}}} M_{\mathfrak{p}} \rbrace \\
& = \sup \lbrace \operatorname{depth} R_{\mathfrak{p}} - \operatorname{width}_{R_{\mathfrak{q}}} M_{\mathfrak{p}} \mid \mathfrak{p} \text{ with } \operatorname{width}_{R_{\mathfrak{p}}} ((B_s)_{\mathfrak{p}}) > \operatorname{width}_{R_{\mathfrak{p}}} M_{\mathfrak{p}} \rbrace \\
& = \sup \lbrace \operatorname{depth} R_{\mathfrak{p}} - \operatorname{width}_{R_{\mathfrak{q}}} M_{\mathfrak{p}} \mid \mathfrak{p} \in \operatorname{Spec} R \rbrace, 
\end{align*}
as required
\end{proof}
\end{theorem}
\begin{remark}
When $C=R$ one recovers from Theorem \ref{Cchouinard} the recent result \cite[Theorem 3.3]{Tri}. Moreover, we remark that the proof of Theorem \ref{Cchouinard} follows a different approach from that of \cite[Theorem 3.3]{Tri}, avoiding results about the small restricted injective dimension.
\end{remark}
We now show that the $C$-quasi-injective dimension of a module coincides with the $C$-injective dimension whenever the latter is finite. The dual statement, concerning projective and quasi-projective dimensions with respect to a semidualizing module, was established for finitely generated modules in \cite[Corollary 5.14]{ferraro}. 

\begin{corollary}\label{crl:comparativo}
Let $C$ be a semidualizing $R$-module and let $M$ be an $R$-module. Then $\operatorname{\textit{C}-qid}_R M \leq \operatorname{\textit{C}-id}_R M$ and equality holds when $\operatorname{\textit{C}-id}_R M$ is finite.
\end{corollary}
\begin{proof}
The inequality $\operatorname{\textit{C}-qid}_RM \leq \operatorname{\textit{C}-id}_R M$ is trivial. Let $\operatorname{\textit{C}-id}_RM<\infty$. Then $M \in \mathcal{A}_C(R)$ by Theorem \ref{corollary2.9}. Therefore, we have  $\Tor_{>0}^R(C,M)=0$. If $\operatorname{\textit{C}-qid}_RM>0$, then using Theorem \ref{Cchouinard} and Proposition \ref{lemma1}, we have
\begin{align*}
\operatorname{\textit{C}-qid}_R M = \sup \lbrace \operatorname{depth} R_{\mathfrak{p}}- \operatorname{width}_{R_{\mathfrak{p}}} M_{\mathfrak{p}} \mid \mathfrak{p} \in \operatorname{Spec}R \rbrace=\operatorname{\textit{C}-id}_RM,
\end{align*}
and the desired equality is proved.

Now, we assume that $\operatorname{\textit{C}-qid}_RM=0$. Since $\operatorname{id}_R (C \otimes_R M) = \operatorname{\textit{C}-id}_RM< \infty$ (using Theorem \ref{theorem2.11}), for each $\mathfrak{p} \in \operatorname{Spec} R$, we have $\operatorname{Ext}_R^{ \gg 0}(R/\mathfrak{p}, C \otimes_R M)=0$. Moreover, as $M \in \mathcal{A}_C(R)$, then $ \operatorname{qid}_R (C \otimes_R M)=\operatorname{\textit{C}-qid}_R M=0$, by \cite[Theorem 4.16]{ferraro}. We must then have that $\operatorname{Ext}_R^{1}(R/\mathfrak{p},C \otimes_R M)=0$ for all $\mathfrak{p} \in \operatorname{Spec} R $, by \cite[Proposition 3.4(2)]{Gheibi}. Thus, $\operatorname{\textit{C}-id}_R M =\operatorname{id}_R (C \otimes_R M) =0$, by \cite[Corollary 3.1.12]{bruns}, as desired.
\end{proof}

\begin{remark}
When $C=R$ one recovers from Corollary \ref{crl:comparativo} the recent result \cite[Proposition 3.6] {jorge2025operatorname}.
\end{remark}
   In the next theorem, we specialize to finitely generated \(R\)-modules. We obtain a more refined result that does not require the assumption \(\operatorname{\textit{C}\text{-}qid}_RM>0\). The theorem holds over an arbitrary ring and extends Bass' formula for the \(C\)-quasi-injective dimension established in \cite[Theorem 7.4]{ferraro} for modules satisfying \(\operatorname{Tor}_{>0}^R(C,M)=0\). Notable this statement and the idea of its proof are inspired from \cite[Theorem 3.1]{Tri}. When \(C=R\), it recovers \cite[Theorem 3.1]{Tri}.
\begin{theorem}\label{theoremfg}
Let $C$ be a semidualizing $R$-module and let $M$ be a finitely generated $R$-module of finite $C$-quasi-injective dimension with $\operatorname{Tor}_{>0}^R(C,M)=0$. Then
$$\operatorname{\textit{C}-qid}_RM = \sup \lbrace \operatorname{depth} R_{\mathfrak{p}} \mid \mathfrak{p} \in \operatorname{Supp}_RM \rbrace.$$
\begin{proof} Let $I_{\bullet}$ be a $C$-quasi-injective resolution of $M$ such that $ \operatorname{\textit{C}-qid}_RM = \operatorname{hinf}(\operatorname{Hom}_R(C,I_{\bullet})) - \operatorname{inf}(\operatorname{Hom}_R(C,I_{\bullet}))$. Without loss of generality, shifting the complex $\operatorname{Hom}_R(C,I_{\bullet})$, we may assume that $\sup(\operatorname{Hom}_R(C, I_{\bullet})) = 0$. Set $-s= \operatorname{hinf}(\operatorname{Hom}_R(C,I_{\bullet}))$. Again, we consider the exact sequences
\begin{align}\label{nonlocal2}
\begin{cases}
0 \rightarrow Z_i \rightarrow \operatorname{Hom}_R(C,I_i) \rightarrow B_{i-1} \rightarrow 0 \\  0 \rightarrow B_i \rightarrow Z_i \rightarrow H_i(\operatorname{Hom}_R(C,I_{\bullet})) \rightarrow 0
\end{cases}
\quad (i \in \mathbb{Z})
\end{align}
where $Z_i=\operatorname{ker}(\partial_i^{\operatorname{Hom}_R(C,I_{\bullet})})$, $B_i=\operatorname{Im}(\partial_{i+1}^{\operatorname{Hom}_R(C,I_{\bullet})})$ and $H_i(\operatorname{Hom}_R(C,I_{\bullet})) \cong M^{\oplus b_i}$ for some $b_i \geq 0$. It is clear that $\operatorname{\textit{C}-qid}_RM=\operatorname{\textit{C}-id}_R Z_{-s}$. As all the injective $R$-modules $I_i$ ($i \in \mathbb{Z}$) are in the Bass class $\mathcal{B}_{C}(R)$ (\cite[1.9]{takahashi}), we must then have $\operatorname{Tor}_{>0}^R(C,\operatorname{Hom}_R(C,I_i))=0$. As $\operatorname{Tor}_{>0}^R(C,M)=0$, then using the exact sequences (\ref{nonlocal2}) one checks by induction that $\operatorname{Tor}_{>0}^R(C,B_i)=0=\operatorname{Tor}_{>0}^R(C,Z_i)$, for all $i$. 

By Proposition \ref{lemma1}, we have
\begin{align}\label{eq:parcial}
    \operatorname{\textit{C}-qid}_R M  = \operatorname{\textit{C}-id}_R Z_{-s}  = \sup \lbrace \operatorname{depth} R_{\mathfrak{p}} - \operatorname{width}_{R_{\mathfrak{q}}} ((Z_{-s})_{\mathfrak{p}}) \mid \mathfrak{p} \in \operatorname{Spec} R \rbrace. 
\end{align}

\par Let $\mathfrak{p} \in \Spec R$. First, we set $\mathfrak{p} \in \Supp_R M$. Consider the short exact sequence $$ 0 \rightarrow B_{-s} \rightarrow Z_{-s} \rightarrow M^{\oplus b_{-s}} \rightarrow 0,$$ where $b_{-s}>0$ since $-s = \operatorname{hinf}(\operatorname{Hom}_R(C,I_{\bullet}))$. This exact sequence localizes into the following exact sequence of $R_{\mathfrak{p}}$-modules
\begin{align*}
      0 \rightarrow (B_{-s})_{\mathfrak{p}} \rightarrow (Z_{-s})_{\mathfrak{p}} \rightarrow M^{\oplus b_{-s}}_{\mathfrak{p}} \rightarrow 0
\end{align*}
that induces the exact sequence   $ k(\mathfrak{p}) \otimes_{R_{\mathfrak{p}}} (Z_{-s})_{\mathfrak{p}} \rightarrow ( k(\mathfrak{p}) \otimes_{R_{\mathfrak{p}}} M_{\mathfrak{p}})^{ \oplus b_{-s}} \rightarrow0$ where $k(\mathfrak{p})$ is the residue field of $R_{\mathfrak{p}}$. Since $0 \neq M_{\mathfrak{p}}$ is finitely generated, it follows by Nakayama's Lemma that $k(\mathfrak{p}) \otimes_{R_{\mathfrak{p}}} M_{\mathfrak{p}} \neq 0$ and so $k(\mathfrak{p}) \otimes_{R_{\mathfrak{p}}} (Z_{-s})_{\mathfrak{p}} \neq 0$. That is, $\operatorname{width}_{R_{\mathfrak{p}}} (Z_{-s})_{\mathfrak{p}}=0$ for all $\mathfrak{p} \in \operatorname{Supp}_R M$.

Now, if $\mathfrak{p} \notin \operatorname{Supp}_R M$, using the exact sequences $$0 \rightarrow B_i \rightarrow Z_i \rightarrow H_i(\operatorname{Hom}_R(C,I_{\bullet})) \rightarrow 0 \quad (i \in \mathbb{Z})
$$ we see that $(B_i)_{\mathfrak{p}} \cong (Z_i)_{\mathfrak{p}}$ for all $i$.  Since  $\operatorname{Tor}_{>0}^R(C,B_i)=\operatorname{Tor}_{>0}^R(C,Z_i)=0$ for all $i$, tensoring by $C \otimes_R -$ the first exact sequence in (\ref{nonlocal2}) we have the exact sequences 
\begin{align}\label{tensor}
    0 \rightarrow C \otimes_R Z_i \rightarrow C \otimes_R \operatorname{Hom}_R(C,I_i) \rightarrow C \otimes_R B_{i-1} \rightarrow0 \quad (i \in \mathbb{Z})
\end{align}
where $C \otimes_R \operatorname{Hom}_R(C,I_i) \cong I_i$ for all $i$ (\cite[1.9]{takahashi}). Let $T$ be any finitely generated $R_{\mathfrak{p}}$-module of finite projective dimension. Localizing the exact sequences (\ref{tensor}) at $\mathfrak{p}$ and applying $\operatorname{Hom}_{R_{\mathfrak{p}}}(T,-)$ we see that 
\begin{align*}
\operatorname{Ext}_{R_{\mathfrak{p}}}^{j+1}(T,C_{\mathfrak{p}} \otimes_{R_{\mathfrak{p}}} (Z_i)_{\mathfrak{p}}) & \cong
   \operatorname{Ext}_{R_{\mathfrak{p}}}^j (T,C_{\mathfrak{p}} \otimes_{R_{\mathfrak{p}}} (B_{ i-1})_{\mathfrak{p}}) \\
   &  \cong \operatorname{Ext}_{R_{\mathfrak{p}}}^j  (T, C_{\mathfrak{p}} \otimes_{R_{\mathfrak{p}}} (Z_{i-1})_{\mathfrak{p}}) \quad \text{for all } i \in \mathbb{Z} \text{ and } j>0. 
\end{align*}
Since we are assuming that $\sup (\operatorname{Hom}_R(C,I_{\bullet}))=0$, we have $Z_0= M^{\oplus b_0}$ for some $b_0 \geq 0$ and using the above isomorphisms we get
\begin{align*}
    \operatorname{Ext}_{R_{\mathfrak{p}}}^j (T,C_{\mathfrak{p}} \otimes_{R_{\mathfrak{p}}} (Z_{-s})_{\mathfrak{p}}) \cong \operatorname{Ext}_{R_{\mathfrak{p}}}^{j+1} (T,C_{\mathfrak{p}} \otimes_{R_{\mathfrak{p}}} (Z_{-s+1})_{\mathfrak{p}}) \cong \cdots \cong \operatorname{Ext}_R^{j+s} (T, C_{\mathfrak{p}} \otimes_{R_{\mathfrak{p}}} (Z_0)_{\mathfrak{p}})=0 
\end{align*}
for all $j>0$ since $(Z_0)_{\mathfrak{p}} = M^{\oplus b_0}_{\mathfrak{p}} =0$, as $\mathfrak{p} \notin \operatorname{Supp}_R M$. As $T$ is any finitely generated $R_{\mathfrak{p}}$-module of finite projective dimension, this shows that $\operatorname{rid}_{R_{\mathfrak{p}}} (C_{\mathfrak{p}} \otimes_{R_{\mathfrak{p}}} (Z_{-s})_{\mathfrak{p}}) \leq 0$ for all $\mathfrak{p} \notin \Supp_R M$. By Fact \ref{factrid}, it implies that $$\operatorname{depth} R_{\mathfrak{p}}- \operatorname{width}_{R_{\mathfrak{p}}} (C_{\mathfrak{p}} \otimes_{R_{\mathfrak{p}}} (Z_{-s})_{\mathfrak{p}}) \leq 0 $$ for all $\mathfrak{p} \notin \Supp_R M$. Moreover, using the same argument that in the proof of Proposition \ref{lemma1}, one can see that $\operatorname{width}_{R_{\mathfrak{p}}} (C_{\mathfrak{p}} \otimes_{R_{\mathfrak{p}}} (Z_{-s})_{\mathfrak{p}})=\operatorname{width}_{R_{\mathfrak{p}}} (Z_{-s})_{\mathfrak{p}} $ and then $\operatorname{depth} R_{\mathfrak{p}}- \operatorname{width}_{R_{\mathfrak{p}}}  (Z_{-s})_{\mathfrak{p}} \leq 0$ for  $\mathfrak{p} \notin \operatorname{Supp}_R M$. Therefore, using the equality (\ref{eq:parcial}), we have 
\begin{align*}
    \operatorname{\textit{C}-qid}_R M & =  \sup \lbrace \operatorname{depth} R_{\mathfrak{p}} - \operatorname{width}_{R_{\mathfrak{q}}} ((Z_{-s})_{\mathfrak{p}}) \mid \mathfrak{p} \in \operatorname{Spec} R \rbrace \\
        & = \sup \lbrace \operatorname{depth} R_{\mathfrak{p}} - \operatorname{width}_{R_{\mathfrak{q}}} ((Z_{-s})_{\mathfrak{p}}) \mid \mathfrak{p} \in \operatorname{Supp}_R M \rbrace
    \\
    & =  \sup \lbrace \operatorname{depth} R_{\mathfrak{p}}  \mid \mathfrak{p} \in \operatorname{Supp}_R M \rbrace.
    \end{align*}
\end{proof}
\end{theorem}

\begin{corollary}
Let $C$ be a semidualizing $R$-module and let $M$ be a finitely generated $R$-module of finite $C$-injective dimension. Then 
$$\operatorname{\textit{C}-id}_RM = \sup \lbrace \operatorname{depth} R_{\mathfrak{p}}\mid \mathfrak{p} \in \operatorname{Supp}_R M \rbrace.$$
\end{corollary}
\begin{proof}
 Since $\operatorname{\textit{C}-id}_RM<\infty$, $M \in \mathcal{A}_C(R)$ by Theorem \ref{corollary2.9}. Then  $\Tor_{>0}^R(C,M)=0$, and the result follows from Corollary \ref{crl:comparativo} and Theorem \ref{theoremfg}.
 \end{proof}
Next, we present a different proof of \cite[Corollary 7.5]{ferraro} that does not make use of Bass' formula and is similar to that of \cite[Corollary 3.2]{Tri}.
\begin{corollary}{\cite[Corollary 7.5]{ferraro}}
Let $(R,\mathfrak{m})$ be a local ring and let $C$ be a semidualizing $R$-module. Let $M$ be a non-zero finitely generated $R$-module such that $\operatorname{\textit{C}-qid}_R M< \infty$ and $\operatorname{Tor}_{>0}^R(C,M)=0$. If $\dim_R M = \dim R$, then $R$ is Cohen-Macaulay.
\begin{proof}
It follows by Theorem \ref{theoremfg} that there exists $\mathfrak{p} \in \operatorname{Supp}_R M$ such that $\operatorname{\textit{C}-qid}_R M= \operatorname{depth} R_{\mathfrak{p}}$. Moreover, \cite[Proposition 7.2]{ferraro} implies that $\operatorname{H}_{\mathfrak{m}}^i(M)=0$ for all $i > \operatorname{\textit{C}-qid}_R M$. Consequently,
\begin{align*}
\dim R =  \dim_R M \leq \operatorname{\textit{C}-qid}_R M = \operatorname{depth} R_{\mathfrak{p}} \leq \operatorname{ht} \mathfrak{p} \leq \dim R,
\end{align*}
where the first inequality  follows from  Grothendieck's Nonvanishing Theorem (\cite[Theorem 3.5.7]{bruns}), since we have $\operatorname{H}_{\mathfrak{m}}^{\dim_R M}(M) \neq 0$.
Hence, $\mathfrak{p}$ must be the maximal ideal of $R$ and therefore $\dim R  \leq \depth R_{\mathfrak{m}}= \depth R$. That is, $R$ is a Cohen-Macaulay ring.
\end{proof}
\end{corollary}
\section{A criterion for finiteness of $C$-injective dimension}
The main theorem of this section yields, as a corollary, a dual version of \cite[Theorem 6.11]{ferraro}, and it recovers the recent result in \cite[Theorem 4.6]{Gheibi} in the case where $C = R$.
\begin{theorem}\label{theorem3}
Let $C$ be a semidualizing $R$-module and let $M$ be an $R$-module such that
\begin{enumerate}
    \item $\operatorname{\textit{C}-qid}_RM< \infty$,
    \item $\operatorname{Tor}_{>0}^R(C,M)=0$,
    \item $\operatorname{Ext}_R^{>0} (C \otimes_RM,C\otimes_RM)=0$ 
\end{enumerate}
then $\operatorname{\textit{C}-id}_RM< \infty$.
\begin{proof}
Let $I_{\bullet}$ be a $C$-quasi-injective resolution of $M$ such that $ \operatorname{\textit{C}-qid}_RM = \operatorname{hinf}(\operatorname{Hom}_R(C,I_{\bullet})) - \operatorname{inf}(\operatorname{Hom}_R(C,I_{\bullet}))$. Without loss of generality, shifting the complex $\operatorname{Hom}_R(C,I_{\bullet})$, we may assume that $\sup(\operatorname{Hom}_R(C, I_{\bullet})) = 0$. Set $s= \operatorname{hinf}(\operatorname{Hom}_R(C,I_{\bullet}))$. As in the proof of Theorem \ref{theoremfg}, we can consider the exact sequences (\ref{nonlocal2}) and check by induction that $\operatorname{Tor}_{>0}^R(C,B_i)=0=\operatorname{Tor}_{>0}^R(C,Z_i)$, as $\Tor_{>0}^R (C,M)=0$. Therefore, by applying $C \otimes_R -$ to the exact sequences (\ref{nonlocal2}) we get exact sequences 
\begin{align}\label{sq6}
\begin{cases}
0 \rightarrow C \otimes_R Z_i \rightarrow C \otimes_R  \operatorname{Hom}_R(C,I_i) \rightarrow C \otimes_R  B_{i-1} \rightarrow 0 \\  0 \rightarrow C \otimes_R  B_i \rightarrow C \otimes_R  Z_i \rightarrow C \otimes_R  H_i(\operatorname{Hom}_R(C,I_{\bullet})) \rightarrow 0
\end{cases}
\quad (i \in \mathbb{Z})
\end{align}
where $C \otimes_R \operatorname{Hom}_R(C,I_i) \cong I_i$ for all $i$ (\cite[1.9]{takahashi}).  It is clear that $\operatorname{\textit{C}-qid}_RM=\operatorname{\textit{C}-id}_R Z_s$ and $\operatorname{\textit{C}-id}_R Z_s=\operatorname{id}_R( C \otimes_R Z_s)<\infty$, by Theorem \ref{theorem2.11}.

By induction, we see that $\operatorname{Ext}_R^{>0}(C \otimes_R M, C \otimes_R Z_i)=0= \operatorname{Ext}_R^{>0}(C \otimes_R M,C\otimes_R B_i)$ for all $i$. Indeed, since we are considering that $\operatorname{sup}(\operatorname{Hom}_R(C,I_{\bullet}))=0$, $Z_0\cong M^{\oplus b_0}$ for some $b_0 \geq 0$. Applying $\operatorname{Hom}_R(C \otimes_R M,-)$ on the exact sequence $0 \rightarrow C \otimes_R Z_0 \rightarrow I_0 \rightarrow C \otimes B_{-1} \rightarrow0$ one can see that $\operatorname{Ext}_R^{>0}(C \otimes_R M, C\otimes_R B_{-1})=0$. Now, considering the exact sequence $$0 \rightarrow C \otimes_R B_{-1} \rightarrow C \otimes_R Z_{-1} \rightarrow C \otimes_R H_{-1}(\operatorname{Hom}_R(C,I_{\bullet})) \rightarrow 0$$ we have $\operatorname{Ext}_R^{>0}(C \otimes_R M,C \otimes_R Z_{-1})=0$. Considering the exact sequences (\ref{sq6}) and repeating this argument we obtain the vanishing of the desired Ext-modules.

Finally, since $H_s(\operatorname{Hom}_R(C,I_{\bullet})) \cong M^{\oplus b_s}$ for some $b_s>0$, we have the short exact sequence
\begin{align*}
0 \rightarrow C \otimes_R B_s \rightarrow C \otimes_R Z_s \rightarrow (C \otimes_R M)^{\oplus b_s} \rightarrow0.
\end{align*}
Since $\operatorname{Ext}_R^1(C \otimes_RM,C \otimes_R B_s)=0$, the above exact sequence splits and $\operatorname{id}_R  (C \otimes_R M)< \infty$, as $\operatorname{id}_R (C \otimes_R Z_{s})< \infty$. Finally, by  Theorem \ref{theorem2.11}, we must then have that $\operatorname{\textit{C}-id}_R M = \operatorname{id}_R (C \otimes_R M)< \infty$.
\end{proof}
\end{theorem}
\begin{remark}
The assumption $\operatorname{Ext}_R^{>0} (C \otimes_RM,C\otimes_RM)=0$ in Theorem \ref{theorem3} can be rewritten as the vanishing of the relative cohomology modules $\operatorname{Ext}_{\mathcal{I}_C}^n(M,M)$ considered in \cite{takahashi} (see \cite[Theorem 4.1]{takahashi}). 
\end{remark}
The following corollary is a dual result to \cite[Theorem 6.11]{ferraro} in the sense of $C$-quasi-injective dimension.
\begin{corollary}\label{crl}
Let $C$ be a semidualizing $R$-module. If $M$ is an $R$-module such that
\begin{enumerate}
    \item $\operatorname{\textit{C}-qid}_R M< \infty$,
    \item $M \in \mathcal{A}_C(R)$,
    \item $\operatorname{Ext}_R^{>0}(M,M)=0$,
\end{enumerate}
then $\operatorname{\textit{C}-id}_R M< \infty$.
\begin{proof}
Since $M \in \mathcal{A}_C(R)$, by \cite[Lemma 3.1.13(a)]{SemidualizingModules}, we have
\begin{align*}
 \operatorname{Ext}_R^i(C\otimes_R M,C \otimes_R M) \cong \operatorname{Ext}_R^i(M,M)=0
\end{align*}
for all $i>0$. By Theorem \ref{theorem3}, it follows that $\operatorname{\textit{C}-id}_RM<\infty$.
\end{proof}
\end{corollary}

\begin{agra}This work was completed while the author was visiting the Department of Mathematics at The University of Texas at Arlington. The author gratefully acknowledges their hospitality. The author thanks Victor D. Mendoza Rubio for helpful discussions regarding this manuscript. The author is grateful to the anonymous referee for a careful reading and for valuable comments.
\end{agra}
\noindent 
\textbf{Conflict of interest.} The corresponding author states that there is no conflict of interest.

\begin{fund}
    The author was supported by grants 2022/12114-0 and 2024/17809-1, S\~ao Paulo Research Foundation (FAPESP).
\end{fund}

\end{document}